\documentclass[11pt,reqno]{amsart}
\usepackage{amssymb,cite}
\numberwithin{equation}{section}

\usepackage[a4paper,
includeheadfoot]{geometry}

\usepackage{latexsym}
\usepackage{textcomp}
\usepackage{amsmath}
\usepackage{amsfonts}
\usepackage{amssymb}
\usepackage{mathrsfs}
\usepackage{enumerate}

\renewcommand{\a }{\alpha }

\renewcommand{\d}{\delta }

\newcommand{\D }{\Delta }

\newcommand{\g }{\gamma}

\renewcommand{\l }{\lambda }

\newcommand{\n }{\nabla }
\newcommand{\vp }{\varphi }
\renewcommand{\phi}{\varphi}

\newcommand{\s }{\sigma }

\renewcommand{\t }{\tau }
\newcommand{\z }{\zeta}

\newcommand{\cO }{\mathcal O}

\newcommand{\be}{\begin{equation}}
\newcommand{\ee}{\end{equation}}

\newcommand{\R}{\mathbb{R}}
\newcommand{\C}{\mathbb{C}}
\newcommand{\N}{\mathbb{N}}
\renewcommand{\S}{\mathbb{S}}

\newcommand{\de}{\partial}

\newcommand{\ti}{\widetilde}

\newcommand{\ra}{{\rangle}}
\newcommand{\la}{{\langle}}

\renewcommand{\k}{\kappa}

\newcommand{\cS}{\mathcal{S}}

\newcommand{\calL }{\mathcal{L}}

\newcommand{\calD }{\mathcal{D}}

\newcommand{\calF}{{\mathcal F}}

\newcommand{\calS}{{\mathcal S}}

\renewcommand{\Re}{{\rm Re}\,}
\renewcommand{\Im}{{\rm Im}\,}

\newtheorem{Theorem}{Theorem}[section]

\newtheorem{Lemma}[Theorem]{Lemma}
\newtheorem{Proposition}[Theorem]{Proposition}
\newtheorem{Corollary}[Theorem]{Corollary}
\newtheorem{Remark}[Theorem]{Remark}

\def\proof{\noindent{{\sc Proof. }}}
\def\square{\vbox{
    \hrule height .4pt
    \hbox{\vrule width .4pt height 7pt \kern 7pt
       \vrule width .4pt}
    \hrule height .4pt }}

\def\square{\vbox{
    \hrule height .4pt
    \hbox{\vrule width .4pt height 7pt \kern 7pt
       \vrule width .4pt}
    \hrule height .4pt }}

\def\QED{\hfill {$\square$}\goodbreak \medskip}

\def\R{{\mathbb R}}

\def\S{{\mathbb S}}

\newcommand{\Ds}{(-\Delta)^s}

\begin{document}
\title {Liouville theorems for a general class of nonlocal operators}
 \author[Mouhamed Moustapha Fall and Tobias Weth]{Mouhamed Moustapha Fall and Tobias Weth}
\address{\hbox{\parbox{5.7in}{\medskip\noindent
      M.M. Fall\\
      African Institute for Mathematical Sciences (A.I.M.S.) of Senegal,
     KM 2, Route de Joal,     B.P. 1418.      Mbour, S\'en\'egal. \\[2pt]
     {\em{E-mail address: }}{\tt mouhamed.m.fall@aims-senegal.org.}\\[5pt]
           }}}
 \address{\hbox{\parbox{5.7in}{\medskip\noindent
      T. Weth\\
       Goethe-Universit\"{a}t Frankfurt, Institut f\"{u}r Mathematik.
Robert-Mayer-Str. 10 D-60054 Frankfurt, Germany. \\[2pt]
     {\em{E-mail address: }}{\tt weth@math.uni-frankfurt.de}\\[5pt]
           }}}
\address{ 
      \textbf{Acknowledgement:} 
M.M.F. is supported by the Alexander von Humboldt Foundation
and partially by the Simons Associateship funding from  the International Center for Theoretical Physics (ICTP). Part of the paper was written while T.W. visited the African Institute for Mathematical Sciences (A.I.M.S.) of Senegal. He wishes to thank the institute for its kind hospitality and the German Academic Exchange Service (DAAD) for funding the visit within the program 57060778.  
}


\date{\today}

 \begin{abstract}
   \noindent
   In this paper, we study the equation $\calL u=0$ in $\R^N$, where $\calL$ belongs to a general class of nonlocal linear operators which may be anisotropic and nonsymmetric.  We classify distributional solutions of this equation, thereby extending and generalizing recent Liouville type theorems in the case where $\calL= (-\Delta)^s$, $s \in (0,1)$ is the classical fractional Laplacian. 
 \end{abstract}

\maketitle
\section{Introduction}
\label{sec:introduction}
In the present paper we consider distributional solutions of operator equations of the form 
$\calL u = 0$, where $\calL$ is related to a class of nonlocal operators $\calL_\nu$ acting on functions $\vp \in C_c^\infty(\R^N)$ via the formula  
\begin{equation}
  \label{eq:general-operator}
[\calL_\nu\vp] (x)= \int_{\R^N}(\vp(x)-\vp(x+y)+y\cdot\n\vp(x)1_{B}(y))\,d\nu(y),\qquad x \in \R^N.
\end{equation}
Here $B$ denotes the unit ball in $\R^N$, and $\nu$ is a signed Radon measure on $\R^N \setminus \{0\}$ with the property that  
\begin{equation}
   \label{eq:general-integral-assumption}
0 < M({\nu}):= \int_{\R^N} \min \{1, |x|^2\} d|\nu|(x) < \infty.
\end{equation}
Here, as usual, $|\nu|$ denotes the associated total variation measure, and integrals over Borel subsets $E \subset \R^N$ with respect to $d\nu$ or $d|\nu|$ will always be understood as integrals over $E \setminus \{0\}$. 
A well studied special case is given by the fractional Laplacian $\calL_\nu= (-\Delta)^s$ with $s \in (0,1)$, which corresponds to the measure 
\be\label{eq:CnS}
d\nu(y)= c_{N,s}|x|^{-N-2s}dy \quad \text{with}\quad  c_{N,s}=  s (1-s)\pi^{-N/2}4^s \frac{\Gamma(\frac{N+2s}{2})}{\Gamma(2-s)}.
\ee
In this case, we may write (\ref{eq:general-operator}) as a principle value integral 
\begin{equation}
  \label{eq:special case}
[\calL_\nu\, \vp] (x)= c_{N,s} PV \int_{\R^N}\frac{\vp(x)-\vp(y)}{|x-y|^{N+2s}}\,dy \qquad 
\text{for $x \in \R^N$,}
\end{equation}
and we have the estimate 
\begin{equation}
  \label{eq:-phi-simple-estimate}
\sup \limits_{x \in \R^N} |[\calL_\nu\, \vp] (x)| (1+|x|^{N+2s})< \infty \qquad \text{for every $\vp \in C_c^\infty(\R^N)$.}
\end{equation}
Consequently, the space 
$$ 
L^1_s(\R^N):= \Bigl \{u \in L^1_{loc}(\R^N)\::\: \int_{\R^N} \frac{|u(x)|}{1+|x|^{N+2s}} \,dx < \infty \Bigr\}
$$
is the natural distributional domain of the operator $(-\Delta)^s$, and we may call a function $u \in L^1_s(\R^N)$ {\em $s$-harmonic in $\R^N$} 
if $\int_{\R^N} u (-\Delta)^s \vp\,dx = 0$ for all $\vp \in C^\infty_c(\R^N)$.  Very recently, it has been shown by the first author in  \cite{F} and independently in 
\cite[Theorem 1.3]{CDL} 
that $s$-harmonic functions in $\R^N$ are affine if $s \in (\frac{1}{2},1)$ and constant 
if $s \in (0,\frac{1}{2}]$. This result generalizes earlier classification theorems stating that bounded or semibounded $s$-harmonic functions are constant, see e.g. \cite{Abat, BKN,FV, RWXZ}. The proof in \cite{F} relies on a Poisson kernel representation 
of  $s$-harmonic functions, whereas the proof of \cite[Theorem 1.3]{CDL} uses Fourier analysis. We point out that the case $s=1$ corresponds to a classical theorem of Joseph Liouville stating that bounded harmonic functions on $\R^N$ are constant. As it is well known, the conclusion also applies to semibounded harmonic functions, see e.g. \cite[Cor. 1.27]{han-lin-pde}. In the present paper, we derive classification results for distributional solutions for a large class of equations involving operators of the type (\ref{eq:general-operator}) which may be anisotropic and nonsymmetric and do not have explicit Poisson kernel representations. As it is common by now, we refer to theorems of this type as Liouville theorems, since they generalize Liouville's classical theorem mentioned above. We will be concerned with  signed Radon measures $\nu$ on $\R^N$ satisfying (\ref{eq:general-integral-assumption}) and, for some $s \ge 0$,  the following decay assumption: 
\begin{equation*}
 |\nu| \bigl(B_1(x)\bigr) = O(|x|^{-N-2s})  \qquad \text{as $|x| \to \infty$.}
 \tag{$D_s$}
\end{equation*}
Note that in the case of absolutely continuous measures given by $d \nu(y)= \kappa(y)dy$ with a function $\kappa \in L^1_{loc}(\R^N)$, assumptions (\ref{eq:general-integral-assumption}) and $(D_s)$ are satisfied if 
\be\label{eq:decay-kappa}
 \int_{B_1(0)} |y|^2 |\kappa(y)|\,dy <\infty \qquad \text{and}\qquad   |\kappa(y)| = O(|y|^{-N-2s}) \qquad \text{as $|x| \to \infty$.} 
\ee 
For a given real number $s \ge 0$, we shall see in Section~\ref{sec:proof-main-result} below that (\ref{eq:general-integral-assumption}) and $(D_s)$ imply the estimate (\ref{eq:-phi-simple-estimate}), and obviously (\ref{eq:-phi-simple-estimate}) remains true if $\nu$ is replaced by the reflected Radon measure $\ti{\nu}$ given by 
\begin{equation}
  \label{eq:definition-refl-measure}
\ti{\nu}(E):=\nu(-E) \qquad \text{for any Borel set $E \subset \R^N \setminus \{0\}$.}
\end{equation}
Moreover, it is easy to see that 
$$
\int_{\R^N} [\calL_{\nu} \psi]   \phi\,dx = \int_{\R^N} \psi \, [\calL_{\ti\nu} \phi]\,dx \qquad \text{for $\psi, \phi \in C^\infty_c(\R^N)$.}
$$
As a consequence, for  $u \in L^1_s(\R^N)$ we may define the distribution $\calL_\nu u$ by 
\begin{equation}
  \label{eq:distr-L}
\langle \calL_\nu u , \phi \rangle := \int_{\R^N} u \, \calL_{\ti\nu} \phi\,dx \qquad \text{for $\phi \in C^\infty_c(\R^N)$.}
\end{equation}
We wish to obtain classification results for distributional solutions of equations containing the operator $\calL_\nu$. Our strategy 
is based on the fact that every $u \in L^1_s(\R^N)$ defines a tempered distribution, so we may characterize $u$ via the distributional support of its Fourier transform $\widehat{u}$. In particular, if $\widehat{u}$ is supported in $\{0\}$, then $u$ is a polynomial (see e.g. \cite[Theorem 6.2]{eskin}). In the special case $\calL_\nu= (-\Delta)^s$, this observation has already been used in \cite{CDL}. At first glance, it is natural to expect that, for a distributional solution $u \in L^1_s(\R^N)$ of $\calL_\nu u = 0$, the 
support of $\widehat{u}$ should be disjoint from the largest open set $\cO \subset \R^N$ where the symbol of $\calL_{\tilde \nu}$ does not vanish. Indeed, this follows easily if the symbol is smooth and $\calL_{\tilde \nu}$ is well defined as an operator on general tempered distributions via Fourier transform. The following result establishes the same property for equations containing $\calL_\nu$ and differential operators under weaker regularity assumptions on the corresponding symbol. 

\begin{Theorem}
\label{main-abstract-result}
Let $s > 0$, let $\nu$ be a signed Radon measure satisfying (\ref{eq:general-integral-assumption}) and $(D_s)$, and let $P$ be a (complex) polynomial. Moreover, let $\cO \subset \R^N$ denote the largest open set such that the symbol 
\begin{equation}
  \label{eq:symbol}
\xi \mapsto \eta(\xi)= -\int_{\R^N}[e^{\imath \xi\cdot y}-1-\imath \xi\cdot y 1_{B(0,1)}(y)]\, d \ti{\nu}(y)
\end{equation}
corresponding to $\ti \nu$ satisfies 
\begin{equation}
  \label{eq:condition-eta}
\eta \in W^{N+2s,1}_{loc}(\cO) \qquad \quad \text{and}\qquad \quad \eta(\xi)+P(-\imath\xi)\neq 0\quad \text{for all $\xi\in \cO$.}  
\end{equation} 
If  $u \in L^1_s(\R^N)$ is a distributional solution of
\be\label{eq:LpP}
\calL_\nu u+P( \n) u=0,
\ee
i.e., 
\be\label{eq:LpP-extended}
\int_{\R^N} u\bigl[\calL_{\ti\nu} \phi +P(- \n) \,\phi\bigr]dx=0 \qquad \text{for all $\phi \in C^\infty_c(\R^N)$,}
\ee
then the  support of  $\widehat{u}$ is contained in  $G:= \R^N \setminus \cO$.\\
Moreover, if $G\subset\{0\}$, then $u$ is a polynomial of degree strictly less than $2s$.\\
\end{Theorem}

\begin{Remark}
\label{remark-on-the-theorem}  
The assumption $\eta \in W^{N+2s,1}_{loc}(\cO)$ implies that, by Sobolev embeddings, $\eta$ is uniquely  represented by a continuous function on $\cO$. Clearly, the second condition in (\ref{eq:condition-eta}) 
 is understood as an assumption on the continuous representation of $\eta$.
\end{Remark}

In the special case where $\nu$ is given by \eqref{eq:CnS} for some $s \in (0,1)$ and thus $\calL_\nu= (-\Delta)^s$ is the fractional Laplacian, the corresponding symbol is given by $\xi \mapsto \eta(\xi)= |\xi|^{2s}$. As a consequence, Theorem~\ref{main-abstract-result} implies the following
result related to L\'evy type operators.

\begin{Corollary}
\label{sec:introduction-1}
Let $s\in (0,1)$, $b\in\R^N$, and let $A \in \R^{N \times N}$ be a positive semidefinite matrix.  If $u\in  L^1_s(\R^N)$ is a distributional solution of 
\begin{equation}
  \label{eq:Levy-type-eq}
\Ds u-{\rm div}(A \nabla u) +b \cdot \n u=0,  
\end{equation}
then there exists $c \in \R$ and $b_* \in \R^N$ such that $b \cdot b_* =0$ and $u(x)= b_* \!\cdot \!x +c$ for $x \in \R^N$. Moreover, $b_*=0$ if $s \le  \frac{1}{2}$. 
\end{Corollary}

To derive this corollary, it suffices to apply Theorem~\ref{main-abstract-result} to $\nu$ given by \eqref{eq:CnS}
and the polynomial $z \mapsto P(z)= - z \cdot Az + b \cdot z$. We then have $\Re (\nu(\xi)+ P(-i\xi))= |\xi|^{2s}+ \xi \cdot A \xi >0$ for $\xi \in \R^N \setminus \{0\}$. 
Thus Theorem~\ref{main-abstract-result}  implies that any distributional solution 
$u\in  L^1_s(\R^N)$ is a polynomial of degree strictly less than $2s$. Hence $u$ is affine, and it is constant if $s \le \frac{1}{2}$. In particular, this implies that $u$ is $s$-harmonic. Writing $u$ in the form $x \mapsto u(x)= b_* \!\cdot \!x +c$, it then follows from (\ref{eq:Levy-type-eq}) that $b \cdot b_* =0$, as claimed.\\

Our main application of Theorem~\ref{main-abstract-result} is concerned with anisotropic variants of the fractional Laplacian. For this we consider the unit sphere $\S^{N-1} \subset \R^N$ and a function $a \in L^\infty(\S^{N-1})$. We then fix $s \in (0,1)$ and let $(-\D)^s_a:=\calL$ be the operator given by (\ref{eq:general-operator}) with 
\begin{equation}
  \label{eq:def-aniso-measure}
d \nu(y)=c_{N,s}  |y|^{-N-2s}a \Bigl(\frac{y}{|y|}\Bigr)\,dy. 
\end{equation}
The assumption $a \in L^\infty(\S^{N-1})$ then ensures that the function 
$\kappa(y):=c_{N,s}  |y|^{-N-2s}a \Bigl(\frac{y}{|y|}\Bigr)$ 
satisfies~\eqref{eq:decay-kappa}
and therefore $\nu$ satisfies $(D_s)$. We have the following result. 
\begin{Theorem}
\label{sec:theorem-anisotropic}
Let $N \ge 2$, and let $a = a_{even} + a_{odd} \in L^\infty(\S^{N-1})$, where $a_{even}$ resp. $a_{odd}$ denote the even and odd part of $a$, respectively. 
Suppose that 
\begin{equation}
  \label{eq:condition}
\int_{\S^{N-1}} | \xi  \cdot \theta  |^{2s}\, a(\theta)\, d \theta >0\qquad \text{for all $\xi \in \S^{N-1}$,}  
\end{equation}
and that 
\begin{equation}
  \label{eq:assumpt-reg-a}
 a_{even} \in W^{\frac{N-1}{2},2}(\S^{N-1}),\qquad   a_{odd} \in W^{\frac{N+2}{2}+2s,2}(\S^{N-1}). 
\end{equation}
Furthermore, let $P$ be a complex polynomial such that $\Re P(- \imath \xi) \ge 0$. Then every distributional solution $u \in L^1_s(\R^N)$ of the equation $(-\D)^s_a u + P(\nabla) u=0$ is affine, and it is constant if $s \le \frac{1}{2}.$
\end{Theorem}

We note that the assumptions of Theorem~\ref{sec:theorem-anisotropic} also include functions which change sign on $\S^{N-1}$. We point out that the regularity assumption in Theorem~\ref{sec:theorem-anisotropic} is weaker in the case where $a \in L^\infty(\S^{N-1})$ is even, and in this case 
the operator $(-\D)^s_a$ is given as a principle value integral 
$$
(-\Delta)^s_a \, \phi (x)= PV \int_{\R^N}\frac{\vp(x)-\vp(y)}{|x-y|^{N+2s}}a\left(\frac{x-y}{|x-y|}\right)dy\qquad \text{for $\phi \in \cS$ and $x \in \R^N$.}
$$ 
Theorem~\ref{sec:theorem-anisotropic} is complementary to a recent interesting Liouville theorem by Ros-Oton and Serra, see \cite[Theorem 2.1]{R-O-S}.  In \cite{R-O-S}, the authors consider anisotropic operators where the function $a \in L^\infty(\S^{N-1})$ above is replaced by an even nonnegative measure on $\S^{N-1}$. In this case, it is in general not possible to define the corresponding operator on the space $L^1_s(\R^N)$ in distributional sense, and instead 
\cite[Theorem 2.1]{R-O-S} relies on the stronger a priori assumption $\|u\|_{L^\infty(B_R(0))} \le C R^\beta$ for $R \ge 1$ with some constants $\beta <2s$ and $C>0$. The argument in \cite{R-O-S} relies on this pointwise growth restriction and does not apply to functions in $L^1_s(\R^N)$.  Our proof of Theorem~\ref{sec:theorem-anisotropic} relies on the well known fact that the real part of the corresponding symbol $\eta$ is homogeneous of degree $2s$, and for $\xi \in \S^{N-1}$ it
is given up to a constant by (\ref{eq:condition}), see Section~\ref{sec:anis-fract-lapl} below. It is an open question whether the regularity assumptions given in (\ref{eq:assumpt-reg-a}) are necessary.

The assumption on $a_{even}$ is related to the fact that, for $\xi \in S^{N-1}$, the expression in (\ref{eq:condition}) is the so-called $2s$-cosine-transformation of $a$. This transformation has nice mapping properties between Hilbertian Sobolev spaces on $\S^{N-1}$ since it is diagonal on spherical harmonics, whereas the corresponding eigenvalues can be computed with the Funk-Hecke formula,  see e.g. \cite{rubin}. In our proof of Theorem~\ref{sec:theorem-anisotropic} in Section~\ref{sec:anis-fract-lapl}, we will also use the Funk-Hecke formula for regularity estimates related to $a_{odd}$ and the imaginary part of the symbol $\eta$.

The paper is organized as follows. Section~\ref{sec:proof-main-result} contains preliminary estimates and the proof of Theorem~\ref{main-abstract-result}. Section~\ref{sec:anis-fract-lapl} is devoted to the family of anisotropic fractional Laplacians and contains the proof of Theorem~\ref{sec:theorem-anisotropic}. In Section~\ref{sec:some-furth-appl}, we briefly present some further applications of Theorem~\ref{main-abstract-result}.

Throughout the paper, we let $\cS$ denote the Schwartz space on $\R^N$ and $\cS'$ the space of tempered distributions. For a tempered distribution $u \in \cS'$, we let both $\widehat u$ and $\calF(u)$ denote the Fourier transform of $u$. Moreover, as usual, $\calF^{-1}(u)$ and $\check u$ stand for the inverse Fourier transform of $u$.

\section{Liouville theorem for L\'{e}vy operators}\label{sec:proof-main-result}
Throughout this section, we assume that $\nu$ is a signed Radon measure on $\R^N$ which satisfies (\ref{eq:general-integral-assumption}) and $(D_s)$ for some $s > 0$. For $k \in \N$, we consider the space  
$$
\calS^k_{s}(\R^N):= \left\{\varphi \in C^k(\R^N) \::\: \sup \limits_{x\in\R^N}(1+|x|^{N+2s}) \sum \limits_{|\alpha| \le k} \:|\partial^{\alpha} \varphi (x)|< \infty\right\} 
$$
endowed with the norm 
$$
\phi \mapsto \|\phi\|_{k,s} := \sup_{x \in \R^N} (1+|x|^{N+2s}) \sum \limits_{|\alpha| \le k} \:|\partial^{\alpha} \varphi (x)|.
$$
\begin{Lemma}
\label{lemma-main-estimate}
Let $\vp \in \calS^2_{s}(\R^N)$.  Then there exists a constant $C>0$ independent of $\vp$ such that \begin{equation}
\int_{\R^N} | \vp(x)-\vp(x+y)-y\cdot\n\vp(x)1_{B(0,1)}(y)|  \,d|\nu|(y) \le C(1+|x|)^{-N-2s} \|\vp\|_{2,s} 
\label{eq:double-inequality}
\end{equation}
for all $x \in \R^N$. Thus $\calL_{\nu} \vp: \R^N \to \R$ is well defined by (\ref{eq:general-operator}) and satisfies 
 \begin{equation}
\label{eq:double-inequality-1}
|\calL_{\nu} \vp(x)| \le  C(1+|x|)^{-N-2s} \|\vp\|_{2,s}\qquad \text{for all $x \in \R^N$.}
\end{equation}
 Moreover, $\calL_{\nu} \vp$ is continuous.
\end{Lemma}

\proof
Let  $\phi \in \cS$. For $x \in \R^N$, we define 
$$
h_x: \R^N \to \R, \qquad h_x(y)= \vp(x)-\vp(x+y)-y\cdot\n\vp(x)1_{B(0,1)}(y).
$$
For $r>0$ and $y \in \R^N$ with $|y| \le r$, we then have, by Taylor expansion, 
\begin{equation}
  \label{eq:r-estimate}
|h_x(y)| \le  \Bigl(\|\n^2\vp\|_{L^\infty(B(x,1))}+2 \|\vp\|_{L^\infty(B(x,r))}\Bigr)   \min \{|y|^2,1\}
\end{equation}
and therefore 
$$
|h_x(y)| \le  3 \|\vp\|_{2,s}  \min \{|y|^2,1\}
\qquad \text{for $x, y \in \R^N$.}
$$
Consequently, 
$$
 \int_{\R^N}|h_x(y)| d|\nu|(y) \le 3 \|\vp\|_{2,s}\, M(\nu)\qquad \text{for $x \in \R^N$,}
$$
and thus, for $x \in \R^N$ with $|x|\le 2$, 
\begin{equation}
  \label{eq:|x|le2}
 \int_{\R^N}|h_x(y)| d|\nu|(y) \le 3^{N+2s+1} \|\vp\|_{2,s}\, M(\nu) (1+|x|)^{-N-2s}.
 \end{equation}
Moreover, for $x \in \R^N$ with $|x|\geq 2$, (\ref{eq:r-estimate}) gives
\begin{align}
\int_{\R^N}&|h_x(y)|\,d|\nu|(y)\leq  \int_{|y|\leq |x|/2}
|h_x(y)|\,d|\nu|(y) + |\vp(x)|\int_{|y| \ge |x|/2}\,d|\nu|(y)\nonumber\\
&\quad \qquad+ \int_{|y|\ge |x|/2}|\vp(x+y)|\,d|\nu|(y)\nonumber\\ 
&\leq  \Bigl(\|\n^2\vp\|_{L^\infty(B(x,1))}+2 \|\vp\|_{L^\infty(B(x,\frac{|x|}{2}))}\Bigr)  M(\nu)
+ |\vp(x)|  M(\nu)\nonumber\\
&\quad \qquad + \|\phi\|_{2,s} \int_{|y| \ge |x|/2}(1+|x+y|)^{-N-2s} \,d|\nu|(y)\nonumber\\
&\leq  \|\vp\|_{2,s}M(\nu) \left( 3\left(1+\frac{|x|}{2}\right)^{-N-2s}  + (1+|x|)^{-
N-2s}\right)  + \|\vp\|_{2,s} H(x) \nonumber\\
&\leq \|\vp\|_{2,s}\Bigl( [3\cdot 2^{N+2s} +1]  M(\nu) (1+|x|)^{-N-2s}   +  H(x)\Bigr),\label{eq:H-1-estimate}
\end{align}
where
$$
H(x):=\int_{|y| \ge 1}(1+|x+y|)^{-
N-2s} \,d|\nu|(y).
$$
To estimate $H(x)$, we use $(D_s)$ and a covering argument. For this we cover $\R^N$ by  disjoint cubes $Q_n$, $n \in \N$ of diameter $1$ (i.e., side length $\frac{1}{\sqrt{N}}$). We assume that one of the cubes, say $Q_0$, is centred at the origin, and we put $\N_*: = \N \setminus \{0\}$.  Then we have the inclusion
$$
\{y \in \R^N\::\: |y| \ge 1\} \subset \bigcup_{n \in \N_*} Q_n.
$$
For $n \in \N_*$ we set 
$r_n:= \inf \limits_{z \in Q_n }|z|$, $R_n:= \sup \limits_{z \in Q_n }|z|$ and 
$$
g_n(x):= \inf_{z \in Q_n}(1+|x+z|)^{-
N-2s},\quad G_n(x):= \sup_{z \in Q_n}(1+|x+z|)^{-
N-2s} \qquad \text{for $x \in \R^N$.}
$$
It is easy to see that there are constants $c_1, c_2>0$ such that $R_n \le c_1 r_n$ for $n \in \N_*$ and 
$$
G_n(x) \le c_2 g_n(x) \qquad \text{for $n \in \N_*$, $x \in \R^N$.}
$$
Moreover, for all $n \in \N$ we have $|Q_n| = N^{-\frac{N}{2}}$, and by (\ref{eq:general-integral-assumption}) and $(D_s$) there exists a constant $c_3>0$ such that $|\nu|(Q_n) \le c_3 r_n^{-N-2s}$ for all $n \in \N_*$. We thus conclude that 
\begin{align*}
\int_{|y| \ge 1}&(1+|x+y|)^{-N-2s} \,d|\nu|(y) \le \sum_{n \in \N^*}\int_{Q_n}(1+|x+y|)^{-N-2s} \,d|\nu|(y)\\
&\le c_3 \sum_{n \in \N_*} G_n(x) r_n^{-N-2s} \le c_1^{N+2s} c_2 c_3 \sum_{n \in \N_*} 
g_n(x) R_n^{-N-2s}\\
&\le c_4 \sum_{n \in \N_*} \int_{Q_ n}(1+|x+y|)^{-N-2s}|y|^{-N-2s}\,dy\\
&\le c_4 \int_{|y| \ge \frac{1}{\sqrt{N}}} (1+|x+y|)^{-N-2s}|y|^{-N-2s}\,dy \qquad \text{with $c_4:=c_1^{N+2s} c_2 c_3$.} 
\end{align*}
Moreover, with 
$$
A_x:= \{y \in \R^N\::\: |y| \ge \frac{1}{\sqrt{N}},\: |x+y| \ge \frac{|x|}{2}\}\;\: \text{and}\;\: 
B_x:= \{y \in \R^N\::\: |y| \ge \frac{1}{\sqrt{N}},\: |x+y| \le \frac{|x|}{2}\},
$$
we find that 
$$
\int_{A_x} (1+|x+y|)^{-N-2s}|y|^{-N-2s}\,dy \le (1+\frac{|x|}{2})^{-N-2s} \int_{|y| \ge \frac{1}{\sqrt{N}} } |y|^{-N-2s}\,dy \le c_5 (1+|x|)^{-N-2s}
$$
and, since $|y| \ge \frac{|x|}{2}$ for $y \in B_x$, 
$$
\int_{B_x} (1+|x+y|)^{-N-2s}|y|^{-N-2s}\,dy \le \bigl(\frac{|x|}{2}\bigr)^{-N-2s} 
\int_{\R^N} (1+|x+y|)^{-N-2s}\,dy= c_6 |x|^{-N-2s}
$$
for $|x| \ge 2$ with $c_5,c_6>0$. Combining these estimates,  we find $c_7>0$ such that $H(x) \le c_7 (1+|x|)^{-N-2s}$ for $x \in \R^N$ with $|x| \ge 2$. Together with (\ref{eq:|x|le2}) and (\ref{eq:H-1-estimate}), it follows that 
$$
\int_{\R^N}|h_x(y)|\,d|\nu|(y) \le C \|\phi\|_{2,s} (1+|x|)^{-N-2s} \qquad \text{for all $x \in \R^N$} 
$$
with a constant $C>0$ independent of $\phi$, as claimed in (\ref{eq:double-inequality}). 
The fact that $\calL_{\tilde \nu}\vp$ is continuous follows from the dominated convergence theorem.
\QED

In the following, we consider the reflected measure $\tilde \nu$ defined by (\ref{eq:definition-refl-measure}), which also satisfies (\ref{eq:general-integral-assumption}) and $(D_s)$. Moreover, we let $\eta$ be the symbol corresponding to $\ti \nu$ as defined in (\ref{eq:symbol}). 

\begin{Lemma}
 \label{simple-fourier-transform}$ $
\begin{itemize}
\item[(i)]  We have $\calF(\calL_{\ti{\nu}} \vp) = \eta \widehat\vp \in C_b(\R^N)$ for every $\vp \in \cS^2_s(\R^N)$.
\item[(ii)] If $P$ is a complex polynomial, $k:= \max \{2, \deg P\}$ and $\vp \in \cS^k_s(\R^N)$ is such that $[\eta   + P (- i \:\cdot)]\widehat\vp \in \cS$, then 
$$
\int_{\R^N} u [\calL_{\ti{\nu}} \vp + P(- \n) \,\phi ] \,dx = \la \widehat u , [\eta   + P (- i \: \cdot)] \widehat\vp \ra  \qquad \text{for every $u \in L^1_s(\R^N)$.}
$$
\end{itemize}
\end{Lemma}

\proof
(i) Let $\vp \in \cS^2_s(\R^N)$. By Lemma~\ref{lemma-main-estimate}, the function 
$$
\R^N \to \R, \qquad x \mapsto  \int_{\R^N} \left| \vp(x)-\vp(x+y)+y\cdot\n\vp(x)1_{B(0,1)}(y)\right|  \,d|\ti{\nu}|(y)
$$
belongs to $L^1(\R^N)$, so that $\calL_{\tilde \nu} \vp \in L^1(\R^N)$  and therefore $\calF(\calL_{\tilde \nu} \vp) \in C_b(\R^N)$. Moreover, Fubini's theorem implies that, for fixed $\xi \in \R^N$, 
\begin{align*}
\calF(\calL_{\ti{\nu}} \vp)(\xi) &= (2\pi)^{-N/2} \int_{\R^N} e^{- \imath x\cdot \xi} \int_{\R^N}[\vp(x)-\vp(x+y)+y\cdot\n\vp(x)1_{B(0,1)}(y)] d\ti{\nu}(y)\,dx\\
&= (2\pi)^{-N/2} \int_{\R^N} \int_{\R^N} e^{- \imath x\cdot \xi}[\vp(x)-\vp(x+y)+y\cdot\n\vp(x)1_{B(0,1)}(y)]dx \,d \ti{\nu}(y)\\
&= \widehat\vp(\xi) \int_{\R^N} [1- e^{\imath y\cdot  \xi} + \imath \,\xi \!\cdot \! y \, 1_{B(0,1)}(y) ]   \,d \ti{\nu}(y)= \widehat\vp(\xi) \eta(\xi).
 \end{align*}
This shows that $\calF(\calL_{\ti{\nu}} \vp) = \eta \widehat\vp \in C_b(\R^N)$.\\
(ii) Let $\vp \in \cS^k_s(\R^N)$ be such that $[\eta   + P (- \imath \: \cdot)]\widehat\vp \in \cS$. Then also 
$\calF^{-1}([\eta   + P (- \imath \: \cdot)]\widehat\vp) \in \cS$, whereas $\calF^{-1}([\eta   + P (- \imath \: \cdot)]\widehat\vp) = \calL_{\ti{\nu}} \vp + P(- \n) \,\phi $ as a tempered distribution as a consequence of (i). It thus follows that $\calL_{\ti{\nu}} \vp + P(- \n) \,\phi \in \cS$, and thus for every $u \in L^1_s(\R^N)$ we have, in distributional sense 
$$
\int_{\R^N} u(x) [\calL_{\ti{\nu}} \vp + P(- \n) \,\phi] (x)\,dx = 
\la u, [\calL_{\ti{\nu}} \vp + P(- \n) \,\phi] \ra = \la \widehat u , [\eta   + P (- \imath \: \cdot)] \widehat\vp \ra,
$$
as claimed.
\QED

\begin{Lemma}
\label{sec:liouv-theor-levy-WN+2s1}
Let $\phi \in W^{N+2s,1}(\R^N)$ be a function with bounded support. Then $\check \phi= \calF^{-1}(\phi) \in \calS_{s}^k(\R^N)$ for every $k \in \N$.   
\end{Lemma}

\proof
Since $\phi$ is a continuous function with compact support, it is clear that $\check \phi \in C^\infty(\R^N)$ and that all derivatives of $\phi$ are bounded on $\R^N$.  In the following, we write $N+2s =m + \alpha$ with $m \in \{N,N+1\}$ and $\alpha \in [0,1)$.  Let $P$ be an arbitrary complex polynomial. Since $\check \phi(x)= \widehat\phi(-x)$, it now suffices to show that 
the functions 
$$
\R^N \to \R, \qquad x \mapsto |x|^\alpha x_i^{m} [P(\nabla)\widehat\phi](x), \qquad i=1,\dots,N
$$
are bounded.  For this we first note that $\psi:= P(\imath \:\cdot) \phi \in W^{m+\alpha,1}(\R^N)$, since $\phi \in W^{m+\alpha,1}(\R^N)$ has bounded support. In particular, for 
$\beta \in \N_0^{N}$ with $|\beta | \le m$ we have $\partial^\beta \psi \in L^1(\R^N)$, and thus the functions   
$$
\R^N \to \R, \qquad  x \mapsto x_i^{m} [P(\nabla)\widehat\phi] = \widehat{\:\imath^m \partial_{i}^m \psi\:}(x) \qquad i=1,\dots,N
$$
are bounded. So the claim follows if $\alpha=0$. If $\alpha \in (0,1)$, we fix $i \in \{1,\dots,N\}$ and consider $\tau= \imath^m\, \partial_{i}^m \psi \in W^{\alpha,1}(\R^N)$. We consider $\t_n\in C^\infty_c(\R^N)$ such that $\t_n\to \t$ in $ W^{\alpha,1}(\R^N)$. We have
$$
\calF(|\cdot|^\a \widehat{\tau_n}(\cdot))= (-\D)^{\frac{\a}{2}}\tau_n
$$
 and thus 
$$
|\xi|^\a \widehat{\tau_n}(\xi)\leq \|(-\D)^{\frac{\a}{2}}\tau_n \|_{L^1(\R^N)}\leq C_{N,\a}\|\t_n\|_{W^{\a,1}(\R^N)}\quad \textrm{ for all }\xi\in\R^N,
$$
where $C_{N,\a}$ is a constant depending only on $N$ and $\a$.
Since also $\t_n\to \tau$ in $ L^1(\R^N)$ as $n \to \infty$, we get 
%
%
$$
\sup_{\xi \in \R^N} |\xi |^\alpha  |\widehat\tau(\xi)| \le C_{N,\a} \|\tau\|_{W^{\alpha,1}(\R^N)}.
$$
Hence the function 
$$
\R^N \to \R, \qquad \xi \mapsto |\xi|^{\alpha} \xi_i^m [P(\nabla)\widehat\phi](\xi)  = |\xi|^\alpha \widehat\tau (\xi)
$$
is bounded, as required.
\QED

\noindent{\sc Proof of Theorem \ref{main-abstract-result}(completed)}.\\
To simplify the notation, we will write $\calL$ instead of $\calL_{\ti{\nu}}$, and we let $k:= \max \{2, \deg P\}$. We first show that 
\be\label{eq:LpP-extended-variant}
\int_{\R^N} u\bigl[\calL \phi +P(- \n) \,\phi\bigr]dx=0 \qquad \text{for all $\phi \in \calS_{s}^k(\R^N)$.}
\ee
Since $C^\infty_c(\R^N)$ is dense in $\calS_{s}^k(\R^N)$, there exists a sequence $(\vp_m)_m$ in $C^\infty_c(\R^N)$ such that 
$$
 \bigl \|\vp -\vp_m \bigr \|_{k,s}  \to 0 \qquad \text{as $m \to \infty$.}
$$
Consequently, by Lemma~\ref{lemma-main-estimate} we have that
$$
\Bigl| \int_{\R^N}u \, \calL (\vp-\vp_m) \,dx \Bigr| \le C \|\vp-\vp_m \|_{2,\lambda} 
\int_{\R^N} \frac{|u(x)|}{(1+|x|)^{N+2s}} dx \to 0 \qquad \text{as $m \to \infty$.} 
$$
Moreover, there exists a constant $c>0$ such that 
$$
\Bigl| \int_{\R^N}u \, P(-\n) (\vp-\vp_m) \,dx \Bigr| \le c \|\vp-\vp_m \|_{k,\lambda} \int_{\R^N} \frac{|u(x)|}{(1+|x|)^{N+2s}} dx \to 0 \qquad \text{as $m \to \infty$.} 
$$
Since $\int_{\R^N}u \, (\calL+P(-\n)) \vp_m \,dx=0$ for all $m \in \N$ by assumption, we obtain 
$$
\int_{\R^N}u \, (\calL+P(-\n)) \vp\,dx = 0,
$$
as claimed.\\
Next, we let $\psi\in C^\infty_c(\cO)$, and we let $K$ be the support of $\psi$. Since $\eta+P(-\imath\: \cdot) \in W^{N+2s,1}_{loc}(\R^N) \subset C(\R^N)$ and $\eta(\xi)+P(-\imath \xi) \not= 0$ in $K$, there exists an open neighborhood $U \subset \subset \cO$ of $K$ with $\inf \limits_{U}|\eta+P(-\imath\: \cdot)|>0$ and $\eta+P(-\imath\: \cdot) \in 
W^{N+2s,1}(U)$. By the chain rule, we then deduce that also $\frac{1}{\eta+P(-\imath\: \cdot)} \in W^{N+2s,1}(U)$. Consequently, we may define $\tau \in W^{N+2s,1}(\R^N)$ by $\tau(\xi) = \frac{\psi(\xi)}{\eta(\xi) + P(-\imath\xi)}$. Since $\tau$ has bounded support in $\R^N$, Lemma~\ref{sec:liouv-theor-levy-WN+2s1} implies that $\phi=\calF^{-1}(\tau) \in \calS_s^k(\R^N)$. Moreover, since $[\eta   + P (- i \: \cdot)]\widehat\vp= \psi \in C^\infty_c(\cO) \subset \cS$, we then have by \eqref{eq:LpP-extended-variant} and Lemma \ref{simple-fourier-transform}(ii):
$$
0= \int_{\R^N}u \, (\calL+P(-\n)\phi \,dx = \la \widehat u , [\eta + P(-\imath \:\cdot)]\widehat\phi \ra  = \la 
\widehat u , \psi \ra
$$
Since this holds for every $\psi\in C^\infty_c(\cO)$, the distributional support of $\widehat u$ is a subset of $ G= \R^N \setminus \cO$. In particular if $G \subset \{0\}$, then  $\widehat u$ is a linear combination of derivatives of the Dirac $\delta$-distribution (see e.g. \cite[Theorem 6.2]{eskin}), so that $u$ is a polynomial. Since $u \in L^1_s(\R^N)$, it follows that the degree of $u$ is strictly less than $2s$.
\QED

\section{The anisotropic fractional Laplacian}
\label{sec:anis-fract-lapl}
The present section is devoted to the proof of Theorem~\ref{sec:theorem-anisotropic}.  In the following, we let $N \ge 2$, and we fix $a \in L^\infty(\S^{N-1})$ and $s \in (0,1)$. We let $a_{even}$ resp $a_{odd}$ denote the even and odd part of $a$, respectively. Moreover, we  let $\nu$ be the signed Radon measure defined by (\ref{eq:def-aniso-measure}), and we let $\eta$ be the symbol corresponding to $\ti \nu$ as given by (\ref{eq:symbol}).  We also recall the definition of the constant $c_{N,s}$ in \eqref{eq:CnS}. We need the following regularity properties of the symbol. 

\begin{Proposition}
\label{sec:anis-fract-lapl-1} 
For $\xi \in \R^N$ we have  
\begin{equation*}
\Re \eta(\xi)= \frac{c_{N,s}}{2c_{1,s}}\int_{\S^{N-1}}| \xi  \cdot \theta |^{2s}\, a(\theta)\, d\theta= \frac{c_{N,s}}{2c_{1,s}}\int_{\S^{N-1}}| \xi  \cdot \theta |^{2s}\, a_{even}(\theta)\, d\theta
\end{equation*}
and 
\begin{equation*}
\Im \eta(\xi)=  c_{N,s} \int_{\S^{N-1}} a(\theta) h_s(\xi \cdot \theta) d \theta = c_{N,s} \int_{\S^{N-1}} a_{odd}(\theta) h_s(\xi \cdot \theta) d \theta
\end{equation*}
with 
\begin{equation}
  \label{eq:def-hs}
h_s \in C^\infty(\R), \qquad  h_s(t)=  \int_0^1 \frac{\sin (r t )  -  r t}{r^{1+2s}}\,dr + \int_1^\infty \frac{\sin (r t)}{r^{1+2s}}\,dr.
\end{equation}
Moreover, we have: 
\begin{itemize}
\item[(i)] If $a_{even} \in W^{\tau,2}(\S^{N-1})$ for some $\tau \ge 0$, then $\Re \eta \in W^{\tau+ \frac{N+1}{2}+ 2s,2}_{loc}(\R^{N} \setminus \{0\}) $.\\  
\item[(ii)] If $a_{odd} \in W^{\tau,2}(\S^{N-1})$ for some $\tau  \ge 0$, then $\Im \eta \in W^{\tau+\frac{N-2}{2},2}_{loc}(\R^{N} \setminus \{0\}) $. 
\end{itemize}
\end{Proposition}

\proof
Clearly we have $\eta(0)=0$. Moreover, since $\tilde \nu$ is given by $d \tilde \nu(y)= c_{N,s}  |y|^{-N-2s}a \Bigl(-\frac{y}{|y|}\Bigr)dy$, 
we have, for fixed $\xi \in \R^N \setminus \{0\}$, 
\begin{align*}
- \frac{\eta(\xi)}{c_{N,s}}&= \int_{\R^N}[e^{\imath \xi \cdot y}-1-\imath \xi\!\cdot\! y 1_{B(0,1)}(y)]\, 
d \tilde \nu(y)\\ &= 
\int_{\S^{N-1}} a(-\theta)  
\int_0^\infty  [e^{\imath r \xi \cdot  \theta}-1- \imath r \xi \!\cdot\! \theta 1_{[0,1]}(r)]r^{-1-2s} \, dr d\theta    
\end{align*}
and therefore 
\begin{align*}
\Re \eta(\xi) & = - c_{N,s} \int_{\S^{N-1}}a(-\theta) \int_0^\infty  [\cos (r \xi \theta) - 1] 
r^{-1-2s} \, dr d\theta\\  
&=- c_{N,s} \int_{\S^{N-1}}a(-\theta) |\xi \!\cdot \! \theta|^{2s}\,d\theta 
\int_0^\infty [\cos t - 1]t^{-1-2s}\,dt\\
& = \frac{c_{N,s}}{2c_{1,s}}  
\int_{\S^{N-1}}a(\theta) |\xi \!\cdot \! \theta|^{2s}\,d\theta = \frac{c_{N,s}}{2c_{1,s}}  
\int_{\S^{N-1}}a_{even}(\theta) |\xi \!\cdot \! \theta|^{2s}\,d\theta.
\end{align*}
Moreover,
\begin{align*}
 \frac{\Im \eta(\xi)}{c_{N,s}}&= 
-\int_{\S^{N-1}}a(-\theta) \int_0^\infty  [\sin (r \xi \theta)  -  r \xi \theta 1_{[0,1]}(t)]r^{-1-2s}\,dr\,d\theta\\
&=  \int_{\S^{N-1}}a(\theta) h_s(\xi \cdot \theta) \,d\theta = \int_{\S^{N-1}}a_{odd}(\theta) h_s(\xi \cdot \theta) \,d\theta
\end{align*}
with $h_s$ as in (\ref{eq:def-hs}), as claimed. Here the last two equalities follow from the oddness of the function $h_s$. Moreover, a standard argument based on Lebesgue's theorem shows that $h_s \in C^\infty(\R)$.\\
To prove (i), we assume that $a_{even} \in W^{\tau,2}(\S^{N-1})$ for some $\tau \ge 0$. In this case, the restriction of $\Re \eta$ to $\S^{N-1}$ coincides with the so-called 
$2s$-cosine-transformation of $a_{even}$, and this transformation has nice mapping properties between Sobolev spaces on $\S^{N-1}$, see e.g. \cite{rubin}. In particular, it follows from \cite[Theorem 1.1]{rubin}, applied with $\alpha=2s+1$,  that $\Re \eta|_{\text{\tiny $\S^{N-1}$}} \in W^{\tau+ \frac{N+1}{2}+ 2s,2}(\S^{N-1})$.  Since $\Re \eta$ is homogeneous of degree $2s$, this easily implies that $\Re \eta \in W^{\tau+ \frac{N+1}{2}+ 2s,2}_{loc}(\R^{N} \setminus \{0\})$.\\
To prove (ii),  we first consider a fixed continuous function $h:[-1,1] \to \R$, and we recall that, for a spherical harmonic $Y_l:\S^{N-1}\to\R$ of degree $l \in \{0,1,2,\dots\}$, 
the Funk-Hecke formula (see e.g. \cite[p. 247]{emot:1953}) yields that 
\begin{equation}
\label{eq1app}
 \int_{\S^{N-1}}h(\xi\cdot \theta) Y_l(\theta)\:d\theta=  \mu(l,N)\, Y_l(\xi)
  \qquad \text{for $\xi \in \S^{N-1}$}
\end{equation}
with 
$$
 \mu(l,N) =\kappa_1(N) 
   \frac{\Gamma(l+1)}{\Gamma(l+N-2)}
   \int_{-1}^{1}h(t)P_l^{\frac{N-2}{2}}(t)d\mu^N (t), 
$$
where 
$$
\kappa_1(N):= 2^{N-2}\pi^{\frac{N-2}{2}}\Gamma\left(\frac{N-2}{2}\right),\qquad  
d\mu^N(t)=(1-t^2)^{\frac{N-1}{2}-1}\:dt
$$ 
and $P_l^\nu$ stands for the Gegenbauer polynomial of order $\nu$ and degree $l$ as defined in 
\cite{szego:1959}. It is not difficult to see that 
$$
\int_{-1}^{1}|P_l^{\frac{N-2}{2}}(t)|^2 d\mu^N (t) = \kappa_2(N) \frac{\Gamma(l+N-2)}{(l +\frac{N-2}{2})\Gamma(l+1)}
$$
with a constant $\kappa_2(N)>0$, see e.g. \cite[eq. (2.4)]{carli}.  More precisely, we have $\kappa_2(N)= \frac{\pi 2^{3-N}}{\Gamma^2(\frac{N-2}{2})}$ for $N \ge 3$,  whereas $\kappa_2(2)>0$ depends on the normalization of zero order Gegenbauer polynomials. Consequently, setting 
$$ 
d_{N,h}:= \Bigl( \int_{-1}^{1}|h(t)|^2 d\mu^N (t)\Bigr)^{1/2},  
$$
 the Cauchy-Schwarz inequality implies that 
$$
\mu(l,N) \le  \kappa_1(N) \sqrt{\kappa_2(N)}d_{N,h} \Bigl(\frac{\Gamma(l+1)}{(l +\frac{N-2}{2}) 
\Gamma(l+N-2)}\Bigr)^{1/2} \le \kappa_1(N) \sqrt{\kappa_2(N)}d_{N,h}\: l^{-\frac{N-2}{2}}
$$
From this we immediately deduce the following regularizing property:
 \begin{itemize}
\item[(R)] If $b \in W^{\tau,2}(\S^{N-1})$, then the function $\xi \mapsto \int_{\S^{N-1}} b(\theta)h(\xi \cdot \theta)\,d\theta$ belongs to $W^{\tau+\frac{N-2}{2},2}(\S^{N-1}).$
\end{itemize}
We now assume that $a_{odd} \in W^{\tau,2}(\S^{N-1})$ for some $\tau  \ge 0$. Note that, in polar coordinates $\rho=|\xi|>0$, $\z= \frac{\xi}{|\xi|} \in \S^{N-1}$, the function  
$M:= \frac{\Im \eta}{c_{N,s}}$ is given by 
$$
M(\rho, \z)= \int_{\S^{N-1}} a_{odd}(\theta) h_s(\rho  \z \cdot \theta)\,d\theta \qquad \text{for $\rho>0$, $ \z \in \S^{N-1}$.}
$$
Since $h_s \in C^\infty(\R)$ and $a_{odd} \in L^1(\S^{N-1})$, a standard argument based on Lebesgue's Theorem yields the existence of 
$$
\partial_{\rho}^k M(\rho, \z) = \int_{\S^{N-1}} a_{odd}(\theta) [ \z \cdot \theta]^{k} 
h_s^{(k)}(\rho  \z \cdot \theta)\,d\theta \qquad \text{for $\rho>0$, $\z \in \S^{N-1}$, $k \in \N$.}
$$
Moreover, for fixed $k \in \N$ and $\rho >0$, we may apply Property (R) above to the function $h \in C([-1,1])$, $h(t)=t^k h_s^{(k)}(\rho t)$ to see that 
$$
\partial_{\rho}^k M(\rho,\cdot) \in W^{\tau+\frac{N-2}{2},2}(\S^{N-1}),
$$
whereas the function $\rho \mapsto \|\partial_{\rho}^k M(\rho,\cdot)\|_{W^{\tau+\frac{N-2}{2},2}(\S^{N-1})}$ is bounded on compact subsets of $(0,\infty)$.
It thus follows that $\Im \eta \in W^{\tau+\frac{N-2}{2},2}_{loc}(\R^N \setminus \{0\})$, as claimed.

\QED
\noindent{\sc Proof of Theorem~\ref{sec:theorem-anisotropic}(completed)}.\\
If $a \in L^\infty(\S^{N-1})$ satisfies assumption (\ref{eq:assumpt-reg-a}), then Proposition~\ref{sec:anis-fract-lapl-1} implies that $\eta \in W^{N+2s,2}_{loc}(\R^{N} \setminus \{0\})
\subset W^{N+2s,1}_{loc}(\R^N \setminus \{0\})$. Moreover, by (\ref{eq:condition}) it follows that $\Re \eta>0$ on $\R^N \setminus \{0\}$, and thus also 
$\Re [\eta+P(-\imath \:\cdot)] >0$ on $\R^N \setminus \{0\}$ by our assumption on the complex polynomial $P$. Hence condition~\eqref{eq:condition-eta} is satisfied with $\cO=\R^N \setminus \{0\}$, and Theorem~\ref{main-abstract-result} then implies that every distributional solution $u \in L^1_s(\R^N)$ of the equation $(-\D)^s_a u + P(\nabla) u=0$ is a polynomial of degree strictly less than $2s$. This implies that every such solution is affine, and it is constant if $s \le \frac{1}{2}.$ The proof is finished. \QED
%
%
%
\section{Some further applications }
\label{sec:some-furth-appl}
An immediate application of Theorem \ref{main-abstract-result} is the following uniqueness result.
\begin{Theorem}\label{th:Dsp1}
Let   $s\in (0,1)$ and assume that $u\in L^1_s(\R^N)$ satisfies $\Ds u+ u=0$ on $\R^N$ in distributional sense. Then $u=0$.
\end{Theorem}

\proof
We consider the operator $\calL_\nu= (-\Delta)^s$ with symbol $\xi \mapsto \eta(\xi)= |\xi|^{2s}$ and the polynomial $P \equiv 1$. Then Theorem~\ref{main-abstract-result} applies with $\cO \supset \R^N \setminus \{0\}$. Hence $u$ is affine, and it is constant if $s \le \frac{1}{2}$. As a consequence, $u$ is $s$-harmonic, which implies that $u= - \Ds u = 0$. 
\QED

Our next result is concerned with the one-dimensional fractional Helmholtz equation $\Ds u- u=0$.

\begin{Theorem}\label{th:1d-HElmotz}
Let   $s\in (0,1)$, and assume that $u\in L^1_s(\R)$ satisfies $\Ds u- u=0$ in $\R$ in distributional sense.  Then $u(x)=c_1\cos(x)+c_2\sin(x)$ for $x \in \R$ with some constants $c_i\in\R$.
\end{Theorem}

\proof
We consider the operator $\calL_\nu= (-\Delta)^s$ with symbol $\xi \mapsto \eta(\xi)= |\xi|^{2s}$ and the polynomial $P \equiv -1$. Then Theorem~\ref{main-abstract-result} applies with $\cO \supset \R \setminus \{0,\pm 1\}$. Consequently, the support of $\widehat u$ is contained in $\{0,\pm 1\}$, which, by the same argument as in \cite[Theorem 6.2]{eskin}, implies that there exists polynomials $p_i$, $i=1,2,3$ such that $u(x)=p_1(x)+ p_2(x)e^{ix} + p_3 e^{-ix}$ for $x \in \R$. Since $u \in L^1_s(\R)$ and $s<1$, it follows that $p_i$ is affine for $i=1,2,3$. Hence we have that $\widehat u = \sum \limits_{k=- 1}^1\Bigl( a_k \delta_k + b_k \delta_k' \Bigr)$ with $a_k,b_k \in \C$, where $\delta_k$ denotes the $\delta$-distribution at the point $k $.  We claim that $b_{1}=b_{-1}= 0$. To see this, we consider $\psi \in C_c^\infty((0,\infty)) \subset C_c^\infty(\R)$ with $\psi(1)=1$.  By Lemma \ref{simple-fourier-transform}(ii) and \eqref{eq:LpP-extended-variant}, applied to $\phi:= {\calF}^{-1}(\psi) \in \cS$, we find that 
$$
0 = \int_{\R^N} u \Bigl[(-\Delta)^s \phi - \phi\Bigr]\,dx = \langle \widehat u, [|\xi|^{2s}-1] \psi  \rangle =b_1 \langle \delta_1', [|\xi|^{2s}-1] \psi \rangle = 2s b_1 \psi(1)= 2s b_1.
$$
Hence $b_1=0$, and similarly we find that $b_{-1}=0$. Consequently, we can write $u$ as 
$$
u(x)= a+b x +c_1 \cos x + c_2 \sin x\qquad \text{for $x \in \R$}
$$
with $a,b,c_1,c_2 \in \R$, whereas $b=0$ if $s \le \frac{1}{2}$.  Since $u$ and the functions $\cos, \sin$ solve the equation $\Ds u = u$ in distributional sense and the function $x \mapsto a + bx$ is $s$-harmonic, it follows that $a= b=0$. The claim thus follows.
\QED

Next we consider the relativistic operator $\calL_{\nu}= (-\D+1)^s-1$, which can be written in the form (\ref{eq:general-operator}) with  $d\nu(y)= c_{N,s}|y|^{-\frac{N+2s}{2}} K_{\frac{N+2s}{2}}(|y|)dy$.  
Here $K_\rho$ is the modified Bessel function of the second kind of order $\rho$, and $c_{N,s}$ is given by \eqref{eq:CnS}, see \cite{FF}. Since $K_{\nu}$ decays exponentially, the measure $\nu$ satisfies (\ref{eq:general-integral-assumption}), and $(D_\s)$ holds for any positive $\s$.
\begin{Theorem}  \label{cor:relat} 
Let  $u\in  L^1_\s(\R^N)$ for some $\s>0$. If
\be\label{eq:relat}
\calL_{\nu} u = ((-\D+1)^s-1) u=0 \quad \textrm{on $\R^N$ in distributional sense,}
\ee
then $u$ is  a harmonic polynomial of degree strictly less than $2\sigma$.
\end{Theorem}

\proof
Since the symbol corresponding to $\calL_{\nu}=\calL_{\tilde \nu}$ is given by $\xi \mapsto (|\xi|^2+1)^s-1$, Theorem~\ref{main-abstract-result} applies with $P \equiv 0$ and $\cO= \R^N \setminus \{0\}$. Consequently, $u$ is a polynomial of degree strictly less than $2\s$. It remains to show that $u$ is harmonic, which follows once we have shown that 
\begin{equation}
  \label{eq:harmonic-fourier}
\langle \widehat{u},|\cdot|^2 \psi \rangle = 0 \qquad \text{for all $\psi \in C^\infty_c(\R^N)$.} 
\end{equation}
For this we consider the function 
$$
h \in C^\infty(\R^N),\qquad h(\xi)= \left\{
  \begin{aligned}
&\frac{|\xi|^2}{(1+|\xi|^2)^s-1},&&\qquad \text{$\xi \not =0$;}\\     
&1/s,&&\qquad \text{$\xi  =0$.}     
  \end{aligned}
\right.
$$
If $\psi \in C^\infty_c(\R^N)$ is given, then $\phi= \calF^{-1}( h \psi)  \in \calS$ satisfies the assumptions of Lemma~\ref{simple-fourier-transform}(ii) with $P\equiv 0$, and this implies that 
$$
0 = \int_{\R^N} u \calL_{\tilde \nu} \phi \,dx = \la \widehat{u}, [(1+|\cdot |^2)^s-1] h \phi \rangle = 
 \langle \widehat{u},|\cdot|^2 \psi \rangle.
$$
Hence (\ref{eq:harmonic-fourier}) holds, and the proof is finished.
\QED

Our final application concerns a nonlocal operator $\calL_\nu$ in $\R$ which appears in 
the intermediate long wave equation from fluid mechanics, see e.g.  \cite{ABS}. 
The operator $\calL_{\nu}$ corresponds to the symbol 
$$
\xi \mapsto  \xi\coth(\pi\xi/2)-2/\pi,\qquad \xi \in \R
$$ 
and can be written in the form (\ref{eq:general-operator}) with $d\nu(y)= \frac{1}{2\pi\sinh^2(y)}dy$, see e.g. \cite[page 6]{R.Frank}.  Since the function $\R \to \R, \; y \mapsto \frac{1}{2\pi\sinh^2(y)}$ decays exponentially at infinity and has a singularity of order $-2$ at $y=0$, the measure $\nu$ satisfies (\ref{eq:general-integral-assumption}), whereas $(D_\s)$ holds for any positive $\s$.

\begin{Theorem}  \label{cor:longw} 
Let $\s>0$, and let $u\in  L^1_\s(\R)$ satisfy $\calL_\nu u=0$ in $\R$ in distributional sense. Then $u$ is a  polynomial of degree strictly less than $2\s$.
\end{Theorem}

\proof
The symbol $\xi\mapsto \xi\coth(\frac{\pi\xi}{2})-\frac{2}{\pi}$ is of class $C^\infty$ and nonzero on $\R \setminus \{0\}$.  To see the  latter, it suffices to note that for $\zeta \in \R \setminus \{0\}$ we have $\tanh \zeta \not = \zeta$  and therefore 
$\coth \zeta \not = \frac{1}{\zeta}$. 
By Theorem   \ref{main-abstract-result}, it thus follows that $u$ is a polynomial of degree strictly less than $2\s$. 
\QED

    \label{References}


\begin{thebibliography}
   \footnotesize

\bibitem{ABS} J. P. Albert, J. L. Bona, J.-C. Saut: Model equations for waves in stratified fluids. Proc. R. Soc. Lond. A 453 (1997), 1233-1260.
\bibitem{Abat} N. Abatangelo: Large s-harmonic functions and boundary blow-up
solutions for the fractional laplacian.
  Discrete Contin. Dyn. Syst. A, V. 35; N. 12 (2015)  5555-5607. 

\bibitem{BKN} K. Bogdan, T. Kulczycki, A. Nowak:
Gradient estimates for harmonic and $q$-harmonic functions of
symmetric stable processes. Illinois J. Math. 46 (2002), no. 2,
541-556.



 \bibitem{CDL}W. Chen, L. D'Ambrosio and  Y. Li: Some Liouville theorems for the fractional Laplacian.   Nonlin. Anal. Volume 121,  2015, Pages 370-381.

\bibitem{carli} L. De Carli: Local $L^p$ inequalities for Gegenbauer polynomials,   in  Topics in classical analysis and applications in honor of Daniel Waterman, 73--87, World Sci. Publ., Hackensack, NJ, (2008).

\bibitem{emot:1953}
Erd\'elyi, A., Magnus, W., Oberhettinger, F., and Tricomi, F.: 
  Higher transcendental functions, Vol. 2, McGraw-Hill, New York 1953 

\bibitem{eskin} 
G. Eskin: Lectures on Linear Partial Differential Equations.   Graduate Studies in Mathematics 123. American Mathematical Society, Providence, Rhode Island, 2011.

\bibitem{FF} M. M.  Fall  and V. Felli:  Unique continuation properties for the relativistic Shr\"{o}rdinger operator with singular potential. Discrete Contin. Dyn. Syst. A, V. 35; N. 12 (2015)   5827-5867.
 
\bibitem{F}M. M. Fall: Entire $s$-harmonic functions are affine.  Proc. Amer. Math. Soc.  In press.
%

\bibitem{FV} F. Ferrari and I. Verbitsky: Radial fractional Laplace operators and Hessian inequalities. J. Differential Equations 253 (2012), no. 1, 244-272. 
\bibitem{R.Frank} R. Frank: On the uniqueness of ground states of non-local equations. J. \'E. D. P.
(2011), Expos\'e no V, 10 p.  


\bibitem{han-lin-pde} Han, Q. and Lin, F.H.: Elliptic partial differential equations. Courant Lecture Notes in Mathematics, 1. New York University, Courant Institute of Mathematical Sciences, New York; American Mathematical Society, Providence, RI, 1997.




%
%
%
%

%

  %

%
\bibitem{R-O-S} X. Ros-Oton and J. Serra: Regularity theory for general stable operators. Preprint 2014. http://arxiv.org/abs/1412.3892. 
%
\bibitem{rubin} B. Rubin: Inversion of fractional integrals related to the spherical Radon transform. J. Funct. Anal. 157 (1998), no. 2, 470--487.

\bibitem{RWXZ} R. Zhuo, W. Chen, X. Cui and Z. Yuan: A Liouville theorem for
the fractional Laplacian. Preprint 2014. arXiv:1401.7402.

\bibitem{szego:1959}
G. Szeg\"{o}, : Orthogonal {P}olynomials. Amer. Math. Soc.,
  Coll. Publ. 23, New York 1959.














    \end{thebibliography}
\end{document}